\documentclass{amsart}
\usepackage{lscape}
\usepackage{amsmath,amssymb}

\begin{document}

\title{Squaring the square with integer linear programming}
\author{Sascha Kurz}
\address{\small Department of Mathematics, Physics, and Computer Science,\\ \small University of Bayreuth, 95440 Bayreuth, Germany\\
  \small Tel.: +49-921-557353, Fax: +49-921-557352, email: sascha.kurz@uni-bayreuth.de}
\maketitle
  
\begin{abstract}
  \noindent
  We consider so-called ``squaring the square''-puzzles where a given square (or rectangle) should be dissected into
  smaller squares. For a specific instance of such problems we demonstrate that a mathematically rigorous solution
  can be quite involved. As an alternative to exhaustive enumeration using tailored algorithms we describe the
  general approach of formulating the problem as an integer linear program.
\end{abstract}

\section{Introduction}
Consider a floor tiler with the task of tiling a rectangular room using square tiles
from a given set only. Questions arising are whether it is possible at all and if so to ask
for minimum cost solutions given a target function like e.g.\ the number of used tiles or the
sum of the actual buying costs for each tile.

Even without this economic interpretation the underlying geometric idea of dissecting a rectangle
into smaller squares was the source for many classical geometric puzzles like e.g.\ the following 
``Problem 173'' in \cite{amusements_in_mathematics}:
\begin{quotation}
  \noindent
  ``For Christmas, Mrs.~Potipher Perkins received a very pretty patchwork quilt constructed of 169~square
  pieces of silk material. The puzzle is to find the smallest number of square portions of which the quilt
  could be composed and show how they might be joined together. Or, to put it the reverse way, divide the
  quilt into as few square portions as possible by merely cutting the stitches.''
\end{quotation}

More formally, the generalized version of this problem can be stated as follows: For a given integer $n$, determine
the minimum number $s(n)$ of squares in a tiling of an $n\times n$-square using $i\times i$-squares with $1\le i\le n-1$,
i.e.\ using only other integral squares.

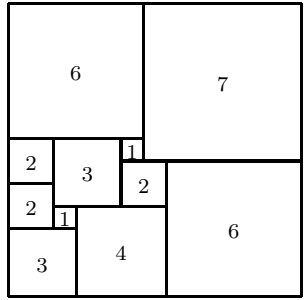
\begin{figure}[htp]
  \begin{center}
    \setlength{\unitlength}{3mm}
      \begin{picture}(13,13)
        \thicklines
        \put(0,0){\line(1,0){13}}   
        \put(0,3){\line(1,0){3}}
        \put(2,4){\line(1,0){5}}
        \put(5,6){\line(1,0){8}}
        \put(0,5){\line(1,0){2}}
        \put(0,7){\line(1,0){6}}
        \put(0,13){\line(1,0){13}}
        \put(0,0){\line(0,1){13}}
        \put(2,3){\line(0,1){4}}
        \put(3,0){\line(0,1){4}}
        \put(5,4){\line(0,1){3}}
        \put(6,6){\line(0,1){7}}
        \put(7,0){\line(0,1){6}}
        \put(13,0){\line(0,1){13}}
        \put(2.25,3.1){\footnotesize{$1$}}
        \put(1.25,1.1){\footnotesize{$3$}}
        \put(0.75,3.6){\footnotesize{$2$}}
        \put(0.75,5.6){\footnotesize{$2$}}
        \put(5.75,4.6){\footnotesize{$2$}}
        \put(4.75,1.6){\footnotesize{$4$}}
        \put(5.25,6.1){\footnotesize{$1$}}
        \put(3.25,5.1){\footnotesize{$3$}}
        \put(2.75,9.6){\footnotesize{$6$}}
        \put(9.25,9.1){\footnotesize{$7$}}
        \put(9.75,2.6){\footnotesize{$6$}}
      \end{picture}
    \caption{Solution of the  patch quilt puzzle.}
    \label{fig_stamp}
  \end{center}
\end{figure}

Due to Martin Gardner this puzzle first appeared in a puzzle magazine edited by Sam Loyd in 1907 and later on
in the famous ``Sam Loyd's Cyclopedia of 5000 Puzzles, Tricks, and Conundrums'', see \texttt{www.mathpuzzle.com/loyd}
for an on-line version. We give the minimal solution consisting of 11~squares in Figure~\ref{fig_stamp} and remark that
it is unique up to rotations and reflections. For an overview on the (scattered) literature concerning these questions
we refer to ``Problem C2'' in \cite{MR1316393}.

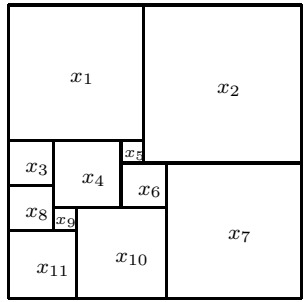
\begin{figure}[htp]
  \begin{center}
    \setlength{\unitlength}{3mm}
      \begin{picture}(13,13)
        \thicklines
        \put(0,0){\line(1,0){13}}   
        \put(0,3){\line(1,0){3}}
        \put(2,4){\line(1,0){5}}
        \put(5,6){\line(1,0){8}}
        \put(0,5){\line(1,0){2}}
        \put(0,7){\line(1,0){6}}
        \put(0,13){\line(1,0){13}}
        \put(0,0){\line(0,1){13}}
        \put(2,3){\line(0,1){4}}
        \put(3,0){\line(0,1){4}}
        \put(5,4){\line(0,1){3}}
        \put(6,6){\line(0,1){7}}
        \put(7,0){\line(0,1){6}}
        \put(13,0){\line(0,1){13}}
        \put(2.1,3.3){\tiny{$x_9$}}
        \put(1.25,1.1){\footnotesize{$x_{11}$}}
        \put(0.75,3.6){\footnotesize{$x_8$}}
        \put(0.75,5.6){\footnotesize{$x_3$}}
        \put(5.75,4.6){\footnotesize{$x_6$}}
        \put(4.75,1.6){\footnotesize{$x_{10}$}}
        \put(5.15,6.25){\tiny{$x_5$}}
        \put(3.25,5.1){\footnotesize{$x_4$}}
        \put(2.75,9.6){\footnotesize{$x_1$}}
        \put(9.25,9.1){\footnotesize{$x_{2}$}}
        \put(9.75,2.6){\footnotesize{$x_{7}$}}
      \end{picture}
    \caption{Unknown sizes of the squares.}
    \label{fig_variables}
  \end{center}
\end{figure}

Once you have found a solution using \textit{few} squares, by a heuristic search or simply by trial and error, it
is easy to verify the validity. Even in the case where the sizes of the squares are omitted one can easily
recover them by solving a linear equation system. Using the variables from Figure~\ref{fig_variables} in
our example this equation system is given by:
\begin{eqnarray*}
  x_1+x_2            &=& 13 \\
  x_3+x_4+x_5+x_2    &=& 13 \\
  x_3+x_4+x_6+x_7    &=& 13 \\
  x_8+x_4+x_6+x_7    &=& 13 \\
  x_8+x_9+x_{10}+x_7 &=& 13 \\
  x_{11}+x_{10}+x_7  &=& 13 \\
  x_1+x_3+x_8+x_{11} &=& 13 \\
  x_1+x_4+x_9+x_{11} &=& 13 \\
  x_1+x_4+x_{10}     &=& 13 \\
  x_1+x_5+x_6+x_{10} &=& 13 \\
  x_2+x_6+x_{10}     &=& 13 \\
  x_2+x_7            &=& 13
\end{eqnarray*}

As in most puzzles, asking for a mini\-mal solution in some sense,
the hardest part is to approve the minimality of the given solution.
This part is addressed by simply stating the smallest \textit{known}
solution and thus not answering this question quite commonly. Only in a few cases
rigorous mathematical proofs are explicitly given and even then circulate
in personal communication only, see e.g.\ 
{\small \texttt{mathworld.wolfram.com/MrsPerkinssQuilt.html}}.

If the puzzle can be formulated as a finite search space, one can in principle
apply exhaustive enumeration. This is the case in the framework of the ``squaring
the square'' context. However, this yields a drawback: sophisticated custom
enumeration algorithms have to be developed and implemented in order to obtain
results in a reasonable amount of time.

The purpose of this article is to demonstrate that it is not that hard to formulate
such puzzles as integer linear programs. Then standard software can
be used in order to exactly solve these problems without a need to implement new algorithms.
The big advantage of such an approach is that it can be easily adapted to different variants
of the puzzle, which we will demonstrate in the following.

\subsection{Outline of the paper}
In Section~\ref{sec_competition_puzzle} we give a puzzle from a ma\-the\-matical
competition for 14--16 year old pupils and outline a rigorous mathematical solution. The known results
on the determination of the minimal number of squares $s(n)$ are outlined in Section~\ref{seq_squaring_the_square}.
The underlying theory for these classes of puzzles is briefly addressed in Section~\ref{sec_electrical_networks}.
In contrast to this exhaustive enumeration approach we describe the modeling process as an integer linear
program in Section~\ref{sec_ilp} and end with a conclusion in Section~\ref{sec_conclusion}.

\section{A puzzle from a mathematical competition and a rigorous solution}
\label{sec_competition_puzzle}
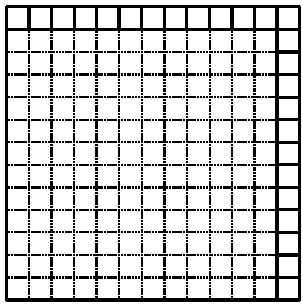
\begin{figure}[htp]
\begin{center}
 \setlength{\unitlength}{3mm}
\begin{picture}(13,13)
\multiput(0,0)(0,1){13}{\qbezier[104](0,0)(6.5,0)(13,0)}
\multiput(0,0)(1,0){13}{\qbezier[104](0,0)(0,6.5)(0,13)}
\thicklines
\put(0,0){\line(1,0){13}}   
\put(0,13){\line(1,0){13}}
\put(0,12){\line(1,0){13}}
\put(12,1){\line(1,0){1}}
\put(12,2){\line(1,0){1}}
\put(12,3){\line(1,0){1}}
\put(12,4){\line(1,0){1}}
\put(12,5){\line(1,0){1}}
\put(12,6){\line(1,0){1}}
\put(12,7){\line(1,0){1}}
\put(12,8){\line(1,0){1}}
\put(12,9){\line(1,0){1}}
\put(12,10){\line(1,0){1}}
\put(12,11){\line(1,0){1}}
\put(0,0){\line(0,1){13}}
\put(13,0){\line(0,1){13}}
\put(12,0){\line(0,1){13}}
\put(1,12){\line(0,1){1}}
\put(2,12){\line(0,1){1}}
\put(3,12){\line(0,1){1}}
\put(4,12){\line(0,1){1}}
\put(5,12){\line(0,1){1}}
\put(6,12){\line(0,1){1}}
\put(7,12){\line(0,1){1}}
\put(8,12){\line(0,1){1}}
\put(9,12){\line(0,1){1}}
\put(10,12){\line(0,1){1}}
\put(11,12){\line(0,1){1}}
\end{picture}
\caption{Optimal solution using a $12\times 12$-square.}
\label{fig_sol_12}
\end{center}
\end{figure}

In a mathematical team competition for 14--16 year pupils we have proposed the following task:
Tile a $13\times 13$ room using $i\times i$-squares where $i=1,2,\dots,12$.
\begin{enumerate}
 \item Determine the minimum number of tiles using at least one $12\times 12$-square.
 \item Determine the minimum number of tiles using at least one $11\times 11$-square.
 \item Determine the minimum number of tiles using at least one $10\times 10$-square.
\end{enumerate}
It is not too hard, this is what we have experienced via the answers of the participants,
to come up with tilings using $26$, $16$, and $13$~squares, respectively. For the first
two questions possible tilings achieving these numbers are drawn in Figures~\ref{fig_sol_12}
and Figure~\ref{fig_sol_11}. Proving that there are no solutions using fewer squares
turned out to be a much harder task for the participants.

To demonstrate the arising difficulties in proving the minimality
of the stated tilings, we outline rigorous proofs for the first two
cases and leave the third case to the interested reader. We suppose
that all squares are arranged on an integer grid.

\subsection{Case 1.}
Using a $12\times 12$-square inside a $13\times13$-square leaves only the possibility
to fill the remaining gaps with $1\times 1$ squares. Since $13^2-12^2=25$ at least $1+25=26$
tiles are required.

\subsection{Case 2.}
\begin{figure}[htp]
\begin{center}
\setlength{\unitlength}{3mm}
\begin{picture}(13,13)
\multiput(0,0)(0,1){13}{\qbezier[104](0,0)(6.5,0)(13,0)}
\multiput(0,0)(1,0){13}{\qbezier[104](0,0)(0,6.5)(0,13)}
\thicklines
\put(0,0){\line(1,0){13}}   
\put(0,13){\line(1,0){13}}
\put(0,11){\line(1,0){11}}
\put(11,2){\line(1,0){2}}
\put(11,4){\line(1,0){2}}
\put(11,6){\line(1,0){2}}
\put(11,8){\line(1,0){2}}
\put(11,10){\line(1,0){2}}
\put(11,11){\line(1,0){2}}
\put(10,12){\line(1,0){1}}
\put(0,0){\line(0,1){13}}
\put(13,0){\line(0,1){13}}
\put(11,0){\line(0,1){11}}
\put(2,11){\line(0,1){2}}
\put(4,11){\line(0,1){2}}
\put(6,11){\line(0,1){2}}
\put(8,11){\line(0,1){2}}
\put(10,11){\line(0,1){2}}
\put(11,11){\line(0,1){2}}
\put(12,10){\line(0,1){1}}
\end{picture}
\caption{Optimal solution using a $11\times 11$-square.}
\label{fig_sol_11}
\end{center}
\end{figure}
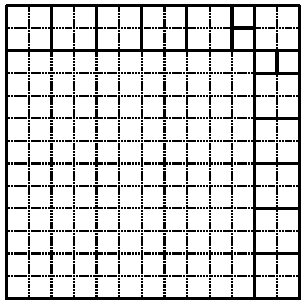

If no corner of the inner $11\times 11$-square would coincide with a corner of the outer
$13\times 13$-square, then at least $13^2-11^2=48$ additional squares are required. Thus
we can assume w.l.o.g.\ that the lower left corners coincide, as in Figure~\ref{fig_sol_11}. As the
largest possible side length of the additional squares is $2$ we can deduce that at least
$\left\lceil\frac{13^2-11^2}{2^2}\right\rceil=12$ of them are needed. Since the
upper $2\times 13$-strip can not be covered using non-overlapping $2\times 2$-squares only
the number $n_2$ of used $2\times 2$-squares is at most $11$ and thus the number $n_1$ of
used $1\times 1$-squares is given by $n_1=13^2-11^2-4n_2\ge 4$. Thus we need at least
$$
  1+n_2+n_1=49-3n_2\ge49-3\cdot11=16 
$$
squares in total.

These ad hoc proofs can be replaced by using a slightly modified version
of the integer linear program presented in Section~\ref{sec_ilp}. We remark that
our proof for the third question is already twice as long as the one stated for case 2.

\section{Squaring the square using as few squares as possible}
\label{seq_squaring_the_square}
Our benchmark example for this article asks for the minimum number $s(n)$ of squares needed
to tile an $n\times n$-square using $i\times i$-squares with $1\le i\le n-1$ only.

Since we can enlarge a given tiling of an $n\times n$-square we obviously have
$s(n_1\cdot n_2)\le \min(s(n_1),s(n_2))$ for all $n_1,n_2\in\mathbb{N}_{\ge 2}$. Thus it
is very likely that we need to determine $s(n)$ for primes only to obtain the minimum values.

In Table~\ref{table_sol_parameters} we list the minimum values $s(p)$ and counts of the used squares
for the first primes. These values were computed by solving the integer linear programing formulation
from Section~\ref{sec_ilp} with the \texttt{Gurobi} solver. Exhaustive enumerations of several kinds
of squared squares are known up to $n=29$, see e.g.\ \texttt{www.squaring.net} and several papers
by A.J.W.\ Duijvestijn like e.g.\ \cite{0502.05017}. Heuristically found configurations, meeting
our exact bounds, can be found in many places. 

The encoding $a_1^{b_1}a_2^{b_2}\dots a_r^{b_r}$, in Table~\ref{table_sol_parameters}, means that
exactly $b_i$ squares of type $a_i\times a_i$ are used.

\setlength{\tabcolsep}{0.2em}
\begin{table}[ht]
  \begin{center}
    \begin{tabular}{ccc}
    \hline
    $\mathbf{p}$ & $\mathbf{s(p)}$ & \textbf{used squares}\\
     $2$ &  $4$ & $1^4$\\
     $3$ &  $6$ & $1^5 2^1$\\
     $5$ &  $8$ & $1^4 2^3 3^1$\\
     $7$ &  $9$ & $1^3 2^3 3^2 4^1$\\
    $11$ & $11$ & $1^4 2^1 3^3 5^2 6^1$\\
    $13$ & $11$ & $1^2 2^3 3^2 4^1 6^2 7^1$\\
    $17$ & $12$ & $1^2 2^3 3^1 4^2 5^1 8^2 9^1$\\
    $19$ & $13$ & $1^1 3^4 4^2 5^3 6^1 9^1 10^1$\\
    $23$ & $13$ & $1^2 2^3 3^1 4^1 6^2 7^1 10^2 13^1$\\
    $29$ & $14$ & $1^1 3^2 4^2 5^3 6^2 7^1 13^2 16^1$\\ 
    $31$ & $15$ & $1^3 2^3 4^3 8^3 15^2 16^1$\\
    $37$ & $15$ & $ 2^2 3^3 4^3 5^1 9^2 11^1 17^2 20^1$\\
    $41$ & $15$ & $1^2 2^2 3^2 4^1 5^1 7^1 11^2 12^1 18^2 23^1$\\
    $43$ & $16$ & $1^2 4^4 5^1 6^2 7^2 11^1 13^1 19^2 24^1$\\
    $47$ & $16$ & $1^1 3^3 6^1 7^2 9^3 9^2 10^1 22^2 25^1$\\
    $53$ & $16$ & $1^1 3^2 4^2 5^2 6^2 7^1 13^2 16^1 24^2 29^1$\\
    $59$ & $17$ & $1^2 2^1 4^2 6^1 8^1 10^2 11^1 12^1 14^1 16^1 18^1 19^1 29^1 30^1$\\
    $61$ & $17$ & $ 3^2 4^4 6^2 8^1 9^2 11^1 15^1 17^1 29^2 32^1$\\
    \hline
    \end{tabular}
    \caption{Minimal numbers $s(p)$ to tile a $p\times p$-square.}
    \label{table_sol_parameters}
  \end{center}
\end{table}

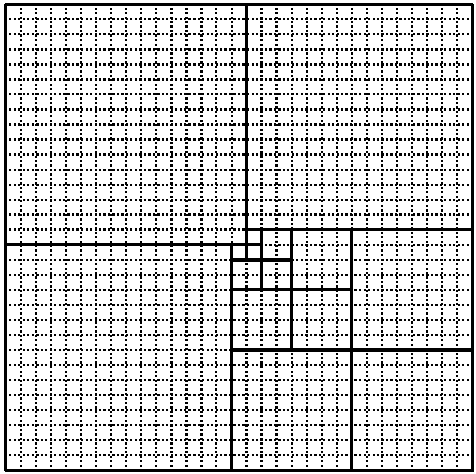
\begin{figure}[ht]
\begin{center}
\setlength{\unitlength}{2mm}
\begin{picture}(31,31)
\multiput(0,0)(0,1){31}{\qbezier[104](0,0)(15.5,0)(31,0)}
\multiput(0,0)(1,0){31}{\qbezier[104](0,0)(0,15.5)(0,31)}
\thicklines
\put(0,0){\line(1,0){31}}   
\put(0,0){\line(0,1){31}}   
\put(31,31){\line(-1,0){31}}
\put(31,31){\line(0,-1){31}}
\put(16,15){\line(-1,0){16}}
\put(16,15){\line(0,1){16}}
\put(15,15){\line(0,-1){15}}
\put(16,16){\line(1,0){15}}
\put(16,15){\line(1,0){1}}
\put(16,15){\line(0,-1){1}}
\put(15,14){\line(1,0){4}}
\put(17,16){\line(0,-1){4}}
\put(15,12){\line(1,0){8}}
\put(19,16){\line(0,-1){8}}
\put(15,8){\line(1,0){16}}
\put(23,16){\line(0,-1){16}}
\end{picture}
\caption{An optimal tiling for a $31\times 31$-square.}
\label{fig_sol_31}
\end{center}
\end{figure}

Looking at $s(p)$ for primes of the form $p=2^r-1$ reveals an interesting pattern:
$1^32^34^3\cdots\left(2^i\right)^3\cdots \left(2^{r-2}\right)^3\left(2^{r-1}-1\right)^2\left(2^{r-1}\right)^1$.\,\, 
In Figure~\ref{fig_sol_31} we have depicted the construction for $p=2^5-1=31$ which can be easily generalized,
so that we have
$$
  s(2^r-1)\le 3r
$$
for all $r\ge 2$.

Trustrum \cite{upper_bound} gave a set of general constructions showing
$$
  s(n)\le 6\log_2(3n-1)-10<6\log_2(n).
$$
For the other direction Conway \cite{lower_bound} has proven
$$
  s(p)\ge \log_2 (p)
$$
for primes $p$.

\section{Squared rectangles from electrical networks}
\label{sec_electrical_networks}
A correspondence between a certain class of planar electrical networks and squared rectangles was observed by
Brooks, Smith, Stone and Tutte \cite{0024.16501}. Here we give a brief sketch of the approach and refer the
interested reader to the expositions \cite{0502.05017,0281.05015} and the historical review \cite{0476.05028}.
Extensive information on the topic can also be found at \texttt{www.squaring.net}.

Given a dissection $\mathcal{D}$ of a rectangle into squares we can build up a network $\mathcal{G}$ as follows. 
Dropping the vertical lines of the constituent squares leaves unions of horizontal lines at the same height, which
we call (horizontal) dissectors. The upper and the lower side of the outer rectangle are examples of such dissectors.
As vertices of $\mathcal{G}$ we choose the horizontal dissectors of the squares and as edges the squares itself. Two
vertices are joined by an edge if the corresponding two horizontal dissectors contain the lower and upper horizontal
side of the corresponding square. We call the vertex corresponding to the upper side of the outer rectangle the positive
pole and the vertex corresponding to the lower side the negative pole of the network.

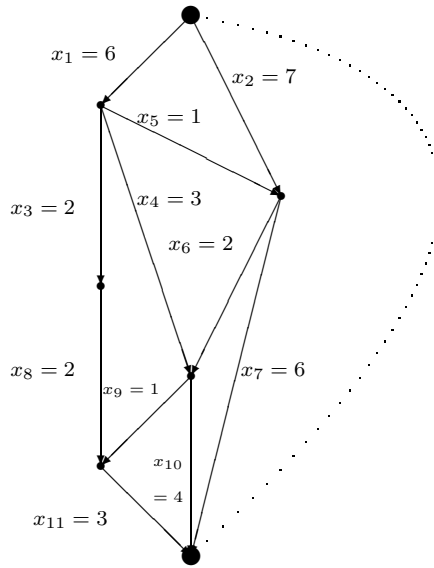
\begin{figure}[htp]
  \begin{center}
    \setlength{\unitlength}{12mm}
    \begin{picture}(2,6)
    \put(1,6){\circle*{0.2}}
    \put(1,6){\vector(1,-2){1}}
    \put(1,6){\vector(-1,-1){1}}
    \put(0,5){\circle*{0.1}}
    \put(0,5){\vector(0,-1){2}}
    \put(0,3){\circle*{0.1}}
    \put(0,3){\vector(0,-1){2}}
    \put(0,1){\circle*{0.1}}
    \put(0,1){\vector(1,-1){1}}
    \put(1,0){\circle*{0.2}}
    \put(2,4){\vector(-1,-2){1}}
    \put(1,2){\circle*{0.1}}
    \put(2,4){\circle*{0.1}}
    \put(0,5){\vector(2,-1){2}}
    \put(0,5){\vector(1,-3){1}}
    \put(1,2){\vector(-1,-1){1}}
    \put(1,2){\vector(0,-2){2}}
    \put(2,4){\vector(-1,-4){1}}
    \put(-0.55,5.5){\footnotesize{$x_1=6$}}
    \put(1.45,5.25){\footnotesize{$x_2=7$}}
    \put(0.4,4.8){\footnotesize{$x_5=1$}}
    \put(-1,3.8){\footnotesize{$x_3=2$}}
    \put(0.4,3.9){\footnotesize{$x_4=3$}}
    \put(0.75,3.4){\footnotesize{$x_6=2$}}
    \put(-1,2.0){\footnotesize{$x_8=2$}}
    \put(0.02,1.8){\tiny{$x_9=1$}}
    \put(-0.75,0.35){\footnotesize{$x_{11}=3$}}
    \put(0.58,1.0){\tiny{$x_{10}$}}
    \put(0.58,0.6){\tiny{$=4$}}
    \put(1.55,2.0){\footnotesize{$x_{7}=6$}}
    \qbezier[60](1,6)(6.5,5)(1,0)
    \end{picture}
    \caption{Network corresponding to Figure~\ref{fig_stamp}.}
    \label{fig_network}
  \end{center}
\end{figure}

The network corresponding to the squared $13\times 13$-square in Figure~\ref{fig_stamp} is drawn in Figure~\ref{fig_network}.

Given a dissection of a rectangle into squares such a network is uniquely defined (if the graph is drawn in a certain way)
and it can be shown that it is planar.
Depending on the properties of the dissection some more graph theoretic restrictions can be deduced. E.g.\ if all squares
have different sizes $\mathcal{G}$ complemented by the edge connecting the positive and the negative pole is a three-connected
planar graph with out multiple edges. If squares of the same size are allowed, the augmented network remains at least two-connected
but may contain multiple edges.

The base for an exhaustive enumeration algorithm for squared rectangles is that from each planar network we can compute the side
lengths of the squares using Kirchhoff's rules. The first law states that the sum of currents in a network meeting at a point is zero.
In our example of Figure~\ref{fig_network} this yields:
\begin{eqnarray*}
  x_1 &=& x_3+x_4+x_5\\
  x_3 &=& x_8\\
  x_2+x_5 &=& x_6+x_{7}\\
  x_4+x_6 &=& x_9+x_{10} \\
  x_8+x_9 &=& x_{11}
\end{eqnarray*}
The second law states that the directed sum of the electrical potential differences in any sub-circuit is zero. Assuming unit resistors
this gives:
\begin{eqnarray*}
  x_1+x_5-x_2 &=& 0 \\
  x_5+x_6-x_4 &=& 0 \\
  x_6+x_{10}-x_{7} &=& 0 \\
  x_9+x_{11}-x_{10}&=&0 \\
  x_3+x_8-x_9-x_4&=&0
\end{eqnarray*}
In general it can be proven that the solution space of the combined equation system is one-dimensional so that one can choose
the unique minimal integer solution. For our example we obtain, of course, multiples of the solution given in Figure~\ref{fig_network}. 

So by exhaustively generating the planar graphs and determining the corresponding dissection of a
rectangle, one can systematically explore the search space for squared rectangles. This approach
is limited to rectangles with a small number of squares. Unless one can exploit strong restrictions on the
graph parameters enumerations for side lengths $n$ with $s(n)\ge 40$ seem to be out of reach,
see \cite{MR2357364}. As already mentioned in Section~\ref{seq_squaring_the_square}, the author is not aware of
any exhaustive enumerations of squared squares with more than $29$ squares. A relatively early work using
computers is \cite{0102.12301}.

\section{An integer linear programming formulation}
\label{sec_ilp}
In this section we want to demonstrate that one can develop an integer linear programming formulation
for the ``squaring the square''-puzzle quite naturally. Once we have such a formulation at hand we can
apply standard software tools to compute the solution. Only minor changes are necessary to adapt the
model to variants of the problem.

To figure out how to model a problem it is useful to ask some basic questions: What can we decide? In our
context the answer is easy -- the positions of the squares. How can we represent or encode our decisions?
Drawing a geometric figure may be suitable for explanations to humans, but talking to a computer we need
something different. We may represent the chosen tiling of Figure~\ref{fig_stamp} by a table:

\smallskip

\begin{center}
     \footnotesize
     \setlength{\tabcolsep}{1mm}
     \begin{tabular}{|c|c|c|c|c|c|c|c|c|c|c|c|c|}
       \hline 6&6&6&6&6&6&7&7&7&7&7&7&7\\[-0.0mm]
       \hline 6&6&6&6&6&6&7&7&7&7&7&7&7\\[-0.0mm]
       \hline 6&6&6&6&6&6&7&7&7&7&7&7&7\\[-0.0mm]
       \hline 6&6&6&6&6&6&7&7&7&7&7&7&7\\[-0.0mm]
       \hline 6&6&6&6&6&6&7&7&7&7&7&7&7\\[-0.0mm]
       \hline 6&6&6&6&6&6&7&7&7&7&7&7&7\\[-0.0mm]
       \hline 2&2&3&3&3&1&7&7&7&7&7&7&7\\[-0.0mm]
       \hline 2&2&3&3&3&2&2&6&6&6&6&6&6\\[-0.0mm]
       \hline 2&2&3&3&3&2&2&6&6&6&6&6&6\\[-0.0mm]
       \hline 2&2&1&4&4&4&4&6&6&6&6&6&6\\[-0.0mm]
       \hline 3&3&3&4&4&4&4&6&6&6&6&6&6\\[-0.0mm]
       \hline 3&3&3&4&4&4&4&6&6&6&6&6&6\\[-0.0mm]
       \hline 3&3&3&4&4&4&4&6&6&6&6&6&6\\[-0.0mm]
       \hline
     \end{tabular}%
\end{center}     

\smallskip

Using this representation the task is to write integers from $1$ to $n-1$ into the cells of a
$n\times n$-grid forming squares for $n=13$. If we have decided where to place the upper left
corner of a $7\times 7$-square, then the remaining $48$ cell entries follow directly. So we
restrict ourselves to printing the positions of the upper left corners:

\smallskip

\begin{center}
\footnotesize
     \setlength{\tabcolsep}{1mm}
 \begin{tabular}{|c|c|c|c|c|c|c|c|c|c|c|c|c|}
       \hline
       6&\phantom{0}&\phantom{0}&\phantom{0}&\phantom{0}&\phantom{0}&7&\phantom{0}&\phantom{0}&\phantom{0}&\phantom{0}&\phantom{0}&\phantom{0}\\[-0.0mm]
       \hline  & & & & & & & & & & & & \\[-0.0mm]
       \hline  & & & & & & & & & & & & \\[-0.0mm]
       \hline  & & & & & & & & & & & & \\[-0.0mm]
       \hline  & & & & & & & & & & & & \\[-0.0mm]
       \hline  & & & & & & & & & & & & \\[-0.0mm]
       \hline 2& &3& & &1& & & & & & & \\[-0.0mm]
       \hline  & & & & &2& &6& & & & & \\[-0.0mm]
       \hline 2& & & & & & & & & & & & \\[-0.0mm]
       \hline  & &1&4& & & & & & & & & \\[-0.0mm]
       \hline 3& & & & & & & & & & & & \\[-0.0mm]
       \hline  & & & & & & & & & & & & \\[-0.0mm]
       \hline  & & & & & & & & & & & & \\[-0.0mm]
       \hline
     \end{tabular}
\end{center}

\smallskip

In many cases binary decision variables are well suited to represent decisions. Therefore we introduce
the binary variables $x_{i,j,h}$ having value $1$ if and only if we write integer $h$ in cell $(i,j)$.

How do we evaluate different decisions? Here our criterion is the number of used squares which can be counted by $$\sum\limits_{i=1}^n\sum\limits_{j=1}^n\sum\limits_{h=1}^{n-1} x_{i,j,h}.$$

What are the constraints restricting our decisions? Here we have to guarantee that each cell of the
$n\times n$-grid is covered exactly by one square and that the squares completely lie inside. The
first condition can be written as
\begin{equation*}
  \sum_{h=1}^{n-1}\sum_{a=0}^{\min(i-1,h-1)}\sum_{b=0}^{\min(j-1,h-1)} x_{i-a,j-b,h}=1
\end{equation*}
for all $1\le i,j\le n$. Understanding this constraint is facilitated by asking the following question: 
Under what condition is a $h\times h$-square with left upper corner at position $(i-a,j-b)$ covering cell $(i,j)$?
The second condition can be written as
$x_{i,j,h}=0$ $\forall 1\le i,j\le n$: $i+h>n+1$ $\vee$ $j+h>n+1$.

Using a modeling language like e.g.\ \texttt{ZIMPL}, see \texttt{zimpl.zib.de}, one can write this in a very compact and readable way:
\begin{verbatim}
param n:=13;
set A:={1 to n};
set B:={1 to n-1};
set S:=A cross A cross B;
var x[S] binary;
minimize target:
sum <i,j,h> in S: x[i,j,h];
subto packing:
forall <i,j> in A cross A:
sum <h> in B:
sum <a> in {0 to min(i-1,h-1)}:
sum <b> in {0 to min(j-1,h-1)}:
x[i-a,j-b,h]==1;
subto boundary:
forall <i,j,h> in S with
i+h>n+1 or j+h>n+1: x[i,j,h]==0;
\end{verbatim}

\texttt{ZIMPL} produces an input file for an ILP solver like \texttt{CPLEX} or \texttt{Gurobi}. In general one can often quite rapidly develop a first working integer linear programming model using \texttt{ZIMPL}, see \cite{Koch2004}. 

\begin{table}[htp]
\begin{center}
     \begin{tabular}{rrrr}
     \hline
      $\mathbf{n}$ & $\mathbf{s(n)}$ & \textbf{seconds} & \textbf{b\&b-nodes}\\
      11 & 11 & 0.1 & 39 \\
      13 & 11 & 0.3 & 41 \\
      17 & 12 & 2   & 92\\
      19 & 13 & 9   & 168\\
      23 & 13 & 19  & 173\\
      29 & 14 & 930 & 3341\\
      30 &  4 & 2   & 1\\
      31 & 15 & 3148& 7409\\
      32 &  4 & 3   & 1\\
      33 &  6 & 4   & 1\\
      34 &  4 & 4   & 1\\
      35 &  8 & 20  & 10\\
      36 &  4 & 5   & 1\\
      37 & 15 & 12634&6911\\
      38 & 4  & 9   & 1\\
      39 & 6  & 8   & 1\\
      40 & 4  & 9   & 1\\
      41 & 15 &26887& 6520\\ 
      \hline
     \end{tabular}
     \caption{Results and running times using \texttt{Gurobi}.}
     \label{table_runtime}
\end{center}
\end{table}

In Table~\ref{table_runtime} we have listed the results for small values of $n$ including the running time and the number of branch\&bound nodes. We have omitted all cases with $n\le 41$, where the problem was solved in the root node, i.e.\ where we have exactly $1$ b\&b-node, and the running time was less than one second. So, without much effort and theoretic insights it was possible to exactly determine the minimum number $s(n)$ of squares to tile a $n\times n$-square, compare sequence A018835 in the on-line encyclopedia of integer sequences at \texttt{oeis.org}.

We remark that using this rather simple approach we could verify the values of Table~\ref{table_sol_parameters} and that $s(n)=\min(s(p)\mid p|n, p\ge 2)$ holds for all $n\le 104$.

\subsection{Problem variations}
Taking the previous integer linear programming formulation as a basic module we can formulate models for variations of the ``squaring the square'' theme. By introducing the auxiliary variables $y_i$ counting the number of used $i\times i$-squares we can express many constraints quite compactly. The relation between the $x$- and the $y$-variables can be stated as
$$
  \sum\limits_{i=1}^n\sum\limits_{j=1}^n x_{i,j,h}=y_h
$$
for all $1\le h\le n-1$.

With this the three problems from Section~\ref{sec_competition_puzzle} can be solved by requiring $y_{12}=1$, $y_{11}=1$, or $y_{10}=1$, respectively.

Since we have observed that the known minimal values $s(n)$ of our benchmark problem from Section~\ref{seq_squaring_the_square} arise from values $s(p)$ for primes $p$ dividing $n$, the additional requirement that the greatest common divisor of the side lengths has to be one was introduced. This can be reformulated such that for each prime $2\le p\le n-1$ the side length of the squares are not all divisible by $p$:
$
  \sum\limits_{1\le h\le n-1\,:\, p\nmid h} y_h \ge 1
$.

\smallskip

\section{Conclusion}
\label{sec_conclusion}
In this article we have considered a special and well known class of geometric puzzles, where rectangles have
to be dissected into smaller squares. Arguing that a discovered solution is minimal in some sense can be quite
involved, as exemplarily demonstrated in Section~\ref{sec_competition_puzzle}. On the other hand performing a
computer based exhaustive search can be time consuming and often requires the development of a customized
algorithm and some theoretical insights in the specific problem. We have briefly outlined the well-known theory
for the problem in Section~\ref{sec_electrical_networks}. Based on exhaustive generation using planar graphs
going beyond $40$ vertices, which correspond to used squares, seems to be computationally infeasible. But we 
have to admit that we are not aware of
any attempts to determine the exact values of $s(n)$ with some kind of restricted generation, i.e.\ where the
search tree is pruned if one can anticipate that the achievable side lengths $n$ will be too large. 

To obtain rigorous results quickly,  modeling the problem as integer linear program and afterwards solving it
with standard software seems to be a viable approach. As demonstrated in Section~\ref{sec_ilp}, the modeling
process is more or less straight-forward and can be adjusted to different problem variations relatively easily.
Of course such a general approach has its computational limits but the same holds for enumeration algorithms too
(whenever the search space grows exponentially, as it does in our context).

Without any sophisticated methods, like e.g.\ column generation, the integer linear programming approach is 
currently not able to push the actual computational frontiers, at least not too much, but seems to be an
accessible way for a broader audience of puzzlers. Its great benefit is its simplicity, compared to the
more involved direct method in Section~\ref{sec_electrical_networks}, and its very general applicability.

Another example where an integer linear programming formulation is used to quickly solve mathematical
puzzles is given in \cite{Koch05}.


\end{document}